# Adaptivity in convolution models with partially known noise distribution

## Cristina Butucea

*Laboratoire Paul Painlevé (UMR CNRS 8524),*
*Université des Sciences et Technologies de Lille 1,*
*59655 Villeneuve d'Ascq cedex, France*
*e-mail:* cristina.butucea@math.univ-lille1.fr

## Catherine Matias

*Laboratoire Statistique et Génome (UMR CNRS 8071),*
*Tour Evry 2, 523 pl. des Terrasses de l'Agora,*
*91000 Evry, France*
*e-mail:* cmatias@genopole.cnrs.fr

## Christophe Pouet

*Laboratoire d'Analyse, Topologie, Probabilités (UMR CNRS 6632),*
*Centre de Mathématiques et Informatique,*
*Université de Provence, 39 rue F. Joliot-Curie,*
*13453 Marseille cedex 13, France*
*e-mail:* pouet@cmi.univ-mrs.fr

**Abstract:** We consider a semiparametric convolution model. We observe random variables having a distribution given by the convolution of some unknown density $f$ and some partially known noise density $g$. In this work, $g$ is assumed exponentially smooth with stable law having unknown self-similarity index $s$. In order to ensure identifiability of the model, we restrict our attention to polynomially smooth, Sobolev-type densities $f$, with smoothness parameter $\beta$. In this context, we first provide a consistent estimation procedure for $s$. This estimator is then plugged-into three different procedures: estimation of the unknown density $f$, of the functional $\int f^2$ and goodness-of-fit test of the hypothesis $H_0 : f = f_0$, where the alternative $H_1$ is expressed with respect to $\mathbb{L}_2$-norm (i.e. has the form $\psi_n^{-2} \|f - f_0\|_2^2 \geq \mathcal{C}$). These procedures are adaptive with respect to both $s$ and $\beta$ and attain the rates which are known optimal for known values of $s$ and $\beta$. These procedures are adaptive with respect to both $s$ and $\beta$ and attain the rates which are known optimal for known values of $s$ and $\beta$. As a by-product, when the noise density is known and exponentially smooth our testing procedure is optimal adaptive for testing Sobolev-type densities. The estimating procedure of $s$ is illustrated on synthetic data.

## 1. Introduction

### *Semiparametric convolution model*

Consider the **semiparametric convolution model** where the observed sample $\{Y_j\}_{1\le j\le n}$ comes from the independent sum of independent and identically distributed (i.i.d.) random variables $X_j$ with unknown density $f$ and Fourier transform $\Phi^f$ and i.i.d. noise variables $\varepsilon_j$ with known, only up to a parameter, density $g$ and Fourier transform $\Phi^g$

$$Y_j = X_j + \varepsilon_j, \quad 1 \le j \le n. \tag{1}$$

The density of the observations is denoted by $p$ and its Fourier transform $\Phi^p$. Note that we have $p = f * g$ where $*$ denotes the convolution product and $\Phi^p = \Phi^f \Phi^g$.

We consider noise distributions whose Fourier transform does not vanish on $\mathbb{R}$: $\Phi^g(u) \ne 0, \forall u \in \mathbb{R}$. Typically, nonparametric estimation in convolution models gives rise to the distinction of two different behaviours for the noise distribution: polynomially or exponentially smooth. In our setup, we focus on exponentially smooth noise where the noise density $g$ may be known only partially. We thus assume an **exponentially smooth** (or supersmooth or exponential) noise distribution having stable symmetric distribution with

$$\Phi^g(u) = \exp\left(-|\gamma u|^s\right), \, \gamma, s > 0. \tag{2}$$

The parameter $s$ is called the self-similarity index of the noise density and we shall consider that it is **unknown** and belongs to a discrete grid $S_n = \{\underline{s} = s_1 < s_2 < \cdots < s_N = \bar{s}\}$, with a number $N$ of points that may grow to infinity with the number $n$ of observations (and $0 < \underline{s} < \bar{s} \le 2$). The parameter $\gamma$ is a scale parameter and it is supposed known in our setting. Some classical examples of such noise densities include the Gaussian and the Cauchy distribution.

The underlying unknown density $f$ is always supposed to belong to $\mathbb{L}_1 \cap \mathbb{L}_2$. For identifiability of the model, the unknown density must be less smooth than the noise. We shall restrict our attention to probability density functions belonging to some Sobolev class

$$\mathcal{S}(\beta, L) = \left\{ f : \mathbb{R} \to \mathbb{R}_+, \int f = 1, \frac{1}{2\pi} \int \left|\Phi^f(u)\right|^2 |u|^{2\beta} du \le L \right\}, \tag{3}$$

for $L$ a positive constant and some unknown $\beta > 0$. We assume that the unknown parameter $\beta$ belongs to some known interval $[\underline{\beta}, \bar{\beta}] \subset (0, +\infty)$. We restrict this interval to $(1/2, +\infty)$ in the case of pointwise estimation of the density $f$. Moreover, we must assume that $f$ is not too smooth, i.e. its Fourier transform does not decay asymptotically faster than a known polynomial of order $\beta'$.

**Assumption (A)** There exists some known $A > 0$, such that $|\Phi^f(u)| \ge A|u|^{-\beta'}$ for large enough $|u|$.

Note that when $f$ belongs to $\mathcal{S}(\beta, L)$ and assumption $(\mathbf{A})$ is fulfilled, we necessarily have $\beta' > \beta + 1/2$. In the following, we use the notation $q_{\beta'}(u) = A|u|^{-\beta'}$. Under assumptions (2) and $(\mathbf{A})$ the model is identifiable. Indeed, considering the Fourier transforms, we get for all real numbers $u$

$$\log |\Phi^p(u)| = \log |\Phi^f(u)| - |\gamma u|^s.$$

Now assume that we have equality between two Fourier transforms of the likelihoods $\Phi_1^p = \Phi_2^p$, where $\Phi_1^p(u) = \Phi^{f_1}(u)e^{-|\gamma_1 u|^{s_1}}$ and $\Phi_2^p(u) = \Phi^{f_2}(u)e^{-|\gamma_2 u|^{s_2}}$. Without loss of generality, we may assume $s_1 \le s_2$. Then we get

$$|u|^{-s_1} \log |\Phi^{f_1}(u)| - \gamma_1^{s_1} = |u|^{-s_2} \log |\Phi^{f_2}(u)| - \gamma_2^{s_2}|u|^{s_2-s_1}$$

and taking the limit when $|u|$ tends to infinity implies (with assumption $(\mathbf{A})$) that $s_1 = s_2$, $\gamma_1 = \gamma_2$ and then $\Phi^{f_1} = \Phi^{f_2}$ which proves the identifiability of the model.

In the sequel, probability and expectation with respect to the distribution of $Y_1, \dots, Y_n$ induced by unknown density $f$ and self-similarity index $s$ will be denoted by $\mathbb{P}_{f,s}$ and $\mathbb{E}_{f,s}$.

Convolution models have been widely studied over the past two decades, mainly in a nonparametric setup where the noise density $g$ is assumed to be entirely known. We will be interested here in a wider framework and will have to deal with the presence of a nuisance parameter $s$. We will focus on both estimation of the unknown density $f$ and goodness-of-fit testing of the hypothesis $H_0 : f = f_0$, with a particular interest in **adaptive** procedures.

Assuming the noise distribution to be entirely known is not realistic in many situations. Thus, dealing with the case of not entirely known noise distribution is a crucial issue. Some approaches [13] rely on additional direct observations from the noise density, which are not always available. A major problem is that semiparametric convolution models do not always result in identifiable models. However, when the noise density is exponentially smooth and the unknown density is restricted to be less smooth than the noise, semiparametric convolution models are identifiable and may be considered.

The case of a Gaussian noise density with unknown variance $\gamma$ and unknown density $f$ without Gaussian component has first been considered in [10]. She proposes an estimator of the parameter $\gamma$ which is then plugged in an estimator of the unknown density. Note that [12] also studied a framework where the variance of the errors is unknown. More generally, [3] consider errors with exponentially smooth stable noise distribution, with unknown scale parameter $\gamma$ but known self-similarity index $s$. The unknown density $f$ belongs either to a Sobolev class, or to a class of supersmooth densities with some parameter $r$ such that $r < s$. Minimax rates of convergence are exhibited. In this context, the unknown parameter $\gamma$ acts as a real nuisance parameter as the rates of convergence for estimating the unknown density are slower compared to the case of known scale, those rates being nonetheless optimal in a minimax sense.

Another attempt to remove knowledge on the noise density appears in [11]. The author proposes a deconvolution estimator associated to a procedure for selecting the error density between the normal supersmooth density and the Laplace polynomially smooth density (both with fixed parameter values). Note that our procedure is more general as we encompass the case of only two different noise distributions and allow a number of unknown supersmooth distributions that may grow to infinity with the number of observations.

Nonparametric goodness-of-fit testing has been extensively studied in the context of direct observations (namely a sample distributed from the density $f$ to be tested), but also for regression or in the Gaussian white noise model. We refer to [9] for an overview on the subject. The convolution model provides an interesting setup where observations may come from a signal observed through some noise.

Nonparametric goodness-of-fit tests in convolution models were studied in [8], [1] and [4], only in the case of entirely known noise distribution. The approach used in [1] is based on a minimax point of view combined with estimation of the quadratic functional $\int f^2$. Assuming the smoothness parameter of $f$ to be known, the authors of [8] define a version of the Bickel-Rosenblatt test statistic and study its asymptotic distribution under the null hypothesis and under fixed and local alternatives, while [1] provides a different goodness-of-fit testing procedure attaining the minimax rates of testing in various setups. The approach used in [1] is further developed in [4] to give adaptive procedures, with respect to the smoothness parameter of $f$, in the case of a polynomially smooth noise distribution.

In our setup, we first propose an estimator of the self-similarity index $s$, which, plugged into kernel procedures, provides an adaptive estimator of the unknown density $f$ with the same optimal rate of convergence as in the case of entirely known noise density. Using the estimator of $s$, we also construct an estimator of the quadratic functional $\int f^2$ (attaining the optimal adaptive rate of convergence) and $\mathbb{L}_2$ goodness-of-fit test statistic. Note that our procedure can only recover the index $s$ on a size-increasing but discrete grid.

Note that this work is very different from [3] as the self similarity index $s$ plays a different role from the scale parameter $\gamma$ previously studied. Nevertheless, we conjecture that their procedure can be extended to recover simultaneously $s$ and $\gamma$ (when both parameters are unknown). However, optimal rates of convergence are even slower when $\gamma$ is unknown.

Another consequence of our results is that when the noise density is known and exponentially smooth our testing procedure is adaptive for testing Sobolev-type densities, improving the previous results in [1].

### *Roadmap*

In Section 2, we provide a consistent estimation procedure for the self-similarity index. Then (Section 3) using a plug-in, we introduce a new kernel estimator of $f$ where both the bandwidth and the kernel are data dependent. We also introduce an estimator of the quadratic functional $\int f^2$ with sample dependent bandwidth and kernel. We prove that these two procedures attain the same rates of convergence as in the case of entirely known noise distribution, and are thus asymptotically optimal in the minimax sense. We also present a goodness-of-fit test on $f$ in this setup. We prove that the testing rate is the same as in the case of entirely known noise distribution and thus asymptotically optimal in the minimax sense. Section 4 illustrates our estimation procedure for parameter $s$ on synthetic data. Proofs are postponed to Section 5.

## 2. Estimation of the self-similarity index $s$

We first present a selection procedure $\hat{s}_n$ which asymptotically recovers the true value of the smoothness parameter $s$ on a given discrete grid

$$S_n = \{\underline{s} = s_1 < s_2 < \cdots < s_N = \bar{s}\},$$

where $0 < \underline{s} < \bar{s} \leq 2$ and with a number $N$ of points that may grow to infinity with $n$, under additional assumptions (see Proposition 1).

Without loss of generality, we assume that $\gamma = 1$ in the following. Indeed, if known $\gamma$ is not equal to 1 then we divide the observations by $\gamma$ to get a noise with scale parameter 1. The asymptotic behavior of the Fourier transform $\Phi^p$ of the observations is used to select the smoothness index $s$. More precisely, we have for any large enough $|u|$

$$A|u|^{-\beta'} \exp(-|u|^s) \leq |\Phi^p(u)| \leq \exp(-|u|^s).$$

Let us now denote $\Phi^{[k]}(u) = e^{-|u|^{s_k}}$ and $I_k(u)$ the interval

$$I_k(u) = [(q_{\beta'}\Phi^{[k]})(u), \Phi^{[k]}(u)],$$

where $q_{\beta'}$ is defined in Assumption (**A**). Let $u_{n,k}$ for $k = 1, \ldots, n$ be some well-chosen points, as described later. Our estimation procedure uses the empirical estimator

$$\hat{\Phi}_n^p(u) = \frac{1}{n} \sum_{j=1}^n \exp(-iuY_j), \quad \forall u \in \mathbb{R},$$

of the Fourier transform $\Phi^p$. We select all values of $k$ belonging to $1, \ldots, N$ such that $\hat{\Phi}_n^p(u_{n,k})$ belongs to or is closest to the interval $I_k(u_{n,k})$. Let then $\hat{s}_n$ be the smallest selected value of $k$, respectively $s_1$ in case no $k$ was selected.

In other words, denote $\hat{S}_n \subset S_n$ the set constructed as follows,

- $s_k \in \hat{S}_n$ if $2 \leq k \leq N - 1$ and

  $$\frac{1}{2} \left\{ q_{\beta'}\Phi^{[k]} + \Phi^{[k+1]} \right\}(u_{n,k}) \leq |\hat{\Phi}_n^p(u_{n,k})| < \frac{1}{2} \left\{ q_{\beta'}\Phi^{[k-1]} + \Phi^{[k]} \right\}(u_{n,k}),$$

- $s_1 \in \hat{S}_n$ if $|\hat{\Phi}_n^p(u_{n,1})| \geq \frac{1}{2} \left\{ q_{\beta'}\Phi^{[1]} + \Phi^{[2]} \right\}(u_{n,1})$,
- $s_N \in \hat{S}_n$ if $|\hat{\Phi}_n^p(u_{n,N})| < \frac{1}{2} \left\{ q_{\beta'}\Phi^{[N-1]} + \Phi^{[N]} \right\}(u_{n,N})$,

where for each index $k$, a sequence of positive real numbers $\{u_{n,k}\}_{n \geq 0}$ has to be chosen later. If the set $\hat{S}_n$ is empty, we add $s_1$. The estimator is

$$\hat{s}_n = \min \hat{S}_n. \tag{4}$$

Note that taking the smallest value such that our condition on the closest interval is satisfied ensures that, with high probability, we do not over-estimate the true value $s$. Over-estimation of $s$ has to be avoided and in some sense, is much worse than under-estimation. Indeed, deconvolution with an over-estimated value of $s$ could result in unbounded estimation risk.

The previous procedure is proven to be consistent, with an exponential rate of convergence, in the following proposition.

**Proposition 1.** *Under assumptions* (2) *and* (**A**), *consider the estimation procedure given by* (4) *where*

$$u_{n,k} = \left( \frac{\log n}{2} - \frac{\delta}{s_k} \log \log n \right)^{1/s_k},$$

*where $\delta > \beta'$. The grid $\underline{s} = s_1 < s_2 < \cdots < s_N = \bar{s}$ is chosen such that*

$$|s_{k+1} - s_k| \geq d_n = \frac{c}{\log n}, \text{ with } c > 2\beta', \quad N - 1 \leq (\bar{s} - \underline{s})/d_n.$$

*Then, for any $k \in \{1, \ldots, N\}$, we have*

$$\mathbb{P}_{f,s_k}(\hat{s}_n \neq s_k) \leq \exp\left( -\frac{A^2}{4} 2^{2\beta'/\bar{s}} (\log n)^{2(\delta - \beta')/\bar{s}} (1 + o(1)) \right), \tag{5}$$

*where $A$ and $\beta'$ are defined in Assumption* (**A**).

## 3. Adaptive estimation and tests

We now plug the preliminary estimator of $s$ in the usual estimation and testing procedures for $f$.

### 3.1. Density estimation

Let us introduce the kernel deconvolution estimator $\hat{K}_n$ (see [5] for a recent survey) built with a preliminary estimation of $s$ plugged-into the usual expression. It is defined by its Fourier transform $\Phi^{\hat{K}_n}$,

$$\Phi^{\hat{K}_n}(u) = \exp\left\{ \left( \frac{|u|}{\hat{h}_n} \right)^{\hat{s}_n} \right\} \mathbb{1}_{|u| \leq 1}, \tag{6}$$

$$\text{where } \hat{h}_n = \left( \frac{\log n}{2} - \frac{\bar{\beta} - \hat{s}_n + 1/2}{\hat{s}_n} \log \log n \right)^{-1/\hat{s}_n}. \tag{7}$$

Note that both the bandwidth sequence $\hat{h}_n$ and the kernel $\hat{K}_n$ are random and depend on observations $Y_1, \ldots, Y_n$. Now, the estimator of $f$ is given by

$$\hat{f}_n(x) = \frac{1}{n\hat{h}_n} \sum_{j=1}^{n} \hat{K}_n \left( \frac{Y_j - x}{\hat{h}_n} \right). \tag{8}$$

This estimation procedure is consistent and adaptively achieves the minimax rate of convergence when considering unknown densities $f$ in the union of Sobolev balls $\mathcal{S}(\beta, L)$ with $\beta \in [\underline{\beta}, \bar{\beta}] \subset (1/2; +\infty)$ and unknown smoothness parameter for the noise density $s$ in a discrete grid $S_n$.

**Corollary 1.** *Under assumptions* (2) *and* (**A**), *for any* $\bar{\beta} > \underline{\beta} > 1/2$, *the estimation procedure given by* (8) *which uses estimator* $\hat{s}_n$ *defined by* (4) *with parameter values:* $\{u_{n,k}\}$ *given by Proposition 1, $\delta > \beta' + \bar{s}^2/(2\underline{s})$, $d_n \geq c \log n$ and $c > 2\beta'$, satisfies, for any real number $x$,*

$$\limsup_{n \to \infty} \sup_{s \in S_n} \sup_{\beta \in [\underline{\beta}, \bar{\beta}]} \sup_{f \in \mathcal{S}(\beta, L)} (\log n)^{(2\beta - 1)/s} \mathbb{E}_{f,s} |\hat{f}_n(x) - f(x)|^2 < \infty.$$

*Moreover, this rate of convergence is asymptotically optimal adaptive.*

**Remark 1.** *This result is obtained by using that, with high probability, the estimator $\hat{s}_n$ is equal to the true value $s$ on the grid (see Proposition 1).*

Note that the optimality of this procedure is a direct consequence of a result by [6] where he considers the convolution model for circular data with $\beta$ and $s$ fixed and known. This result confirms the results of [2] for adaptive estimation of linear functionals in the convolution model and known parameter $s$. Therefore we may say that there is no loss due to adaptation neither with respect to $\beta$ nor to $s$.

Note also that by similar calculations we get that the adaptive estimator $\hat{f}_n$ attains the rate $(\log n)^{2\beta/s}$ over Hölder classes of probability density functions of smoothness $\beta$, for the mean squared error (pointwise risk).

Moreover, it can be shown that the mean integrated squared error of the adaptive estimator $\hat{f}_n$ converges at the rate $(\log n)^{2\beta/s}$ over either Sobolev or Hölder classes of functions. In [7], lower bounds of the same order were proven over Hölder classes of density functions $f$.

### 3.2. Goodness-of-fit test

In the sequel, $\| \cdot \|_2$ denotes the $\mathbb{L}_2$-norm, $\bar{M}$ is the complex conjugate of $M$ and $< M, N > = \int M(x)\bar{N}(x)dx$ is the scalar product of complex-valued functions in $\mathbb{L}_2(\mathbb{R})$. From now on, we consider again that $[\underline{\beta}, \bar{\beta}] \subset (0, +\infty)$.

For a given density $f_0$ in the class $\mathcal{S}(\beta_0, L_0)$, we want to test the hypothesis

$$H_0 : f = f_0$$

from observations $Y_1, \ldots, Y_n$ given by (1). We extend the results of [1] by giving the family of sequences $\Psi_n = \{\psi_{n,\beta}\}_{\beta \in [\underline{\beta}, \bar{\beta}]}$ which separates (with respect to $\mathbb{L}_2$-norm) the null hypothesis from a larger alternative

$$H_1(\mathcal{C}, \Psi_n) : f \in \cup_{\beta \in [\underline{\beta}, \bar{\beta}]} \{f \in \mathcal{S}(\beta, L) \text{ and } \psi_{n,\beta}^{-2} \|f - f_0\|_2^2 \geq \mathcal{C}\}.$$

Let us first remark that as we use noisy observations (and unlike what happens with direct observations), this test cannot be reduced to testing uniformity of the distribution density of the observed sample (i.e. $f_0 = 1$ with support on the finite interval $[0; 1]$).

We recall that the usual procedure is to construct, for any $0 < \epsilon < 1$, a test statistic $\Delta_n^\star$ (an arbitrary function, with values in $\{0, 1\}$, which is measurable with respect to $Y_1, \ldots, Y_n$ and such that we accept $H_0$ if $\Delta_n^\star = 0$ and reject it otherwise) for which there exists some $\mathcal{C}^0 > 0$ such that

$$\limsup_{n \to \infty} \sup_{s \in S_n} \left\{ \mathbb{P}_{f_0, s}[\Delta_n^\star = 1] + \sup_{f \in H_1(\mathcal{C}, \Psi_n)} \mathbb{P}_{f, s}[\Delta_n^\star = 0] \right\} \leq \epsilon, \qquad (9)$$

holds for all $\mathcal{C} > \mathcal{C}^0$. This part is called the upper bound of the testing rate. Then, prove the minimax optimality of this procedure, i.e. the lower bound

$$\liminf_{n \to \infty} \inf_{\Delta_n} \sup_{s \in S_n} \left\{ \mathbb{P}_{f_0, s}[\Delta_n = 1] + \sup_{f \in H_1(\mathcal{C}, \Psi_n)} \mathbb{P}_{f, s}[\Delta_n = 0] \right\} \geq \epsilon, \qquad (10)$$

for some $\mathcal{C}_0 > 0$ and for all $0 < \mathcal{C} < \mathcal{C}_0$, where the infimum is taken over all test statistics $\Delta_n$.

An additional assumption **(T)** used in [1] on the tail behaviour of $f_0$ (ensuring it does not vanish arbitrarily fast) is needed to obtain the optimality result, which is in fact a consequence of [1]. We recall this assumption here for reader's convenience.

**Assumption (T)**

$$\exists c_0 > 0, \forall x \in \mathbb{R}, f_0(x) \geq \frac{c_0}{1 + |x|^2}.$$

**Remark 2.** *Similar results may be obtained under the more general assumption: there exists some $p \geq 1$ such that $f_0(x)$ is bounded from below by $c_0(1 + |x|^p)^{-2}$ for large enough $x$.*

Now, the first step is to construct an estimator of $\int f^2$. Using the same kernel estimator (6) and the same random bandwidth (7), we define

$$\hat{T}_n = \frac{2}{n(n-1)} \sum_{1 \leq k < j \leq n} < \frac{1}{\hat{h}_n} \hat{K}_n \left( \frac{\cdot - Y_k}{\hat{h}_n} \right) , \frac{1}{\hat{h}_n} \hat{K}_n \left( \frac{\cdot - Y_j}{\hat{h}_n} \right) > . \qquad (11)$$

**Corollary 2.** *Under assumptions (2) and (A), for any $\bar{\beta} > \underline{\beta} > 0$, the estimation procedure given by (11) which uses estimator $\hat{s}_n$ defined by (4) with* *parameter values:* $\{u_{n,k}\}$ *given by Proposition 1, $\delta > \beta' + \bar{s}^2/(2\underline{s})$, $d_n \geq c \log n$ and $c > 2\beta'$, satisfies,*

$$\limsup_{n \to \infty} \sup_{s \in \bar{S}_n} \sup_{\beta \in [\underline{\beta}, \bar{\beta}]} \sup_{f \in \mathcal{S}(\beta, L)} (\log n)^{2\beta/s} \left\{ \mathbb{E}_{f,s} \left| \hat{T}_n - \int f^2 \right|^2 \right\}^{1/2} < \infty.$$

*Moreover, this rate of convergence is asymptotically adaptive optimal.*

The rate of convergence of this procedure is the same as in the case of known self-similarity index $s$ and known smoothness parameter $\beta$. It is thus asymptotically optimal adaptive according to results obtained by [1].

Let us now define, for any $f_0 \in \mathcal{S}(\bar{\beta}, L)$,

$$\hat{T}_n^0 = \frac{2}{n(n-1)} \sum_{1 \leq k < j \leq n} < \frac{1}{\hat{h}_n} \hat{K}_n \left( \frac{\cdot - Y_k}{\hat{h}_n} \right) - f_0 , \frac{1}{\hat{h}_n} \hat{K}_n \left( \frac{\cdot - Y_j}{\hat{h}_n} \right) - f_0 > . \tag{12}$$

This statistic is used for goodness-of-fit testing of the hypothesis $H_0$ versus $H_1$. The test is constructed as usual

$$\Delta_n^\star = \begin{cases} 1 & \text{if } |\hat{T}_n^0| \hat{t}_n^{-2} > \mathcal{C}^\star \\ 0 & \text{otherwise,} \end{cases} \tag{13}$$

for some constant $\mathcal{C}^\star > 0$ and a **random** threshold $\hat{t}_n^2$ to be specified.

For computational facilities, we may write using Plancherel formula

$$\hat{T}_n^0$$
$$= \frac{2}{n(n-1)} \sum_{1 \leq k < j \leq n} \frac{1}{2\pi} < \Phi^{\hat{K}_n} \left( \cdot \hat{h}_n \right) e^{i \cdot Y_k} - \Phi^{f_0} , \Phi^{\hat{K}_n} \left( \cdot \hat{h}_n \right) e^{-i \cdot Y_j} - \overline{\Phi^{f_0}} >$$
$$= \frac{1}{\pi n(n-1)} \sum_{1 \leq k < j \leq n}$$
$$\int \left( e^{|u|^{\hat{s}_n} + iuY_k} I_{\hat{h}_n|u| \leq 1} - \Phi^{f_0}(u) \right) \cdot \left( e^{|u|^{\hat{s}_n} - iuY_j} I_{\hat{h}_n|u| \leq 1} - \overline{\Phi^{f_0}}(u) \right) du.$$

**Corollary 3.** *Under assumptions (2) and (A), for any $0 < \underline{\beta} < \bar{\beta}$, any $L > 0$ and for any $f_0 \in \mathcal{S}(\bar{\beta}, L)$, consider the testing procedure given by (13) which uses the test statistic (12) with estimator $\hat{s}_n$ defined by (4) with parameter values: $\{u_{n,k}\}$ given by Proposition 1, $\delta > \beta' + \bar{s}^2/(2\underline{s})$, $d_n \geq c \log n$ and $c > 2\beta'$, with random threshold and (slightly modified) random bandwidth*

$$\hat{t}_n^2 = \left( \frac{\log n}{2} \right)^{-2\bar{\beta}/\hat{s}_n} \quad ; \quad \hat{h}_n = \left( \frac{\log n}{2} - \frac{2\bar{\beta}}{\hat{s}_n} \log \log n \right)^{-1/\hat{s}_n}$$

*and any large enough positive constant $\mathcal{C}^\star$. This testing procedure satisfies (9) for any $\epsilon \in (0, 1)$ with testing rate*

$$\Psi_n = \{\psi_{n,\beta}\}_{\beta \in [\underline{\beta}, \bar{\beta}]} \text{ given by } \psi_{n,\beta} = \left( \frac{\log n}{2} \right)^{-\beta/s}.$$

*Moreover, if $f_0 \in \mathcal{S}(\bar{\beta}, cL)$ for some $0 < c < 1$ and if Assumption* (**T**) *holds, then this testing rate is asymptotically adaptive optimal over the family of classes $\{\mathcal{S}(\beta, L), \beta \in [\underline{\beta}; \bar{\beta}]\}$ and for any $s \in S_n$ (i.e.* (10) *holds).*

Adaptive optimality (namely (10)) of this testing procedure directly follows from [1] as there is no loss due to adaptation to $\beta$ nor to $s$. Note also that the case of known $s$ and adaptation only with respect to $\beta$ is included in our results and is entirely new.

## 4. Simulations

In this section, we illustrate some of our results on synthetic data. We consider two different signal densities: the density of the sum of 5 independent Laplace random variables, *Laplace*(5) (having standard deviation $\sqrt{10}$) and a Gamma distribution with parameters $(3/2, 1/2)$ or $\chi_3^2$ (with standard deviation $\sqrt{6}$), as described in Table 1.

The noise densities were selected among 4 different exponentially smooth distributions as described in Table 2.

The simulation of random variables having Fourier transform $\Phi^{[0.5]}(u)$ and $\Phi^{[1.5]}(u)$ is based on [14]. We thus simulated 8 different samples each one containing $n$ observations, where $n$ ranges from $\{500; 1000; 2000; 5000\}$. We used a scale 0.1 on the signal density in order to have a small signal-to-noise ratio (defined as the ratio of the standard deviations of the signal with respect to that of the noise). Note that the noise has finite standard deviation only for $s = 2$ and it equals $\sqrt{2}$. In this case, the signal-to-noise ratio equals 0.22 when the signal has Laplace density and 0.17 for the Gamma distribution.

We then performed selection of $s$ on the finite grid $S_n = \{0.5, 1, 1.5, 2\}$. The points $u_{n,k}$ were chosen independently of the size $n$ of the sample. The choice is based both on theoretical grounds and on a previous simulation study. We fixed the following values $u_{n,1} = 2.5; u_{n,2} = 1.7; u_{n,3} = 1.5; u_{n,4} = 1.45$. For each sample and each sample size, we performed $m = 100$ iteration of the procedure and the results are presented in Table 3 for the Laplace signal density and in 4 for the Gamma signal density.

TABLE 1
*Signal densities*

| Signal density | Fourier transform |
|---|---|
| *Laplace*(5) | $\Phi_L(u) = (1 + u^2)^{-5}$ |
| *Gamma*$(\frac{3}{2}, \frac{1}{2})$ | $\Phi_G(u) = (1 - 2iu)^{-3/2}$ |

TABLE 2
*Noise densities*

| Noise stable density | Fourier transform |
|---|---|
| $s = 0.5$ | $\Phi^{[0.5]}(u) = \exp(-|u|^{1/2})$ |
| $s = 1$ | $\Phi^{[1]}(u) = \exp(-|u|)$ |
| $s = 1.5$ | $\Phi^{[1.5]}(u) = \exp(-|u|^{1.5})$ |
| $s = 2$ | $\Phi^{[2]}(u) = \exp(-|u|^2)$ |

TABLE 3
*Number of successes ($\hat{s}_n = s$) for 100 iterations of the procedure, when the signal density is Laplace*

|         | $n = 500$ | $n = 1000$ | $n = 2000$ | $n = 5000$ |
|---------|-----------|------------|------------|------------|
| $s = 0.5$ | 85 | 93 | 98 | 100 |
| $s = 1$   | 66 | 87 | 95 | 100 |
| $s = 1.5$ | 65 | 82 | 93 | 100 |
| $s = 2$   | 73 | 90 | 93 | 99 |

TABLE 4
*Number of successes ($\hat{s}_n = s$) for 100 iterations of the procedure, when the signal density is Gamma*

|         | $n = 500$ | $n = 1000$ | $n = 2000$ | $n = 5000$ |
|---------|-----------|------------|------------|------------|
| $s = 0.5$ | 94 | 99 | 100 | 100 |
| $s = 1$   | 71 | 88 | 98 | 100 |
| $s = 1.5$ | 91 | 98 | 100 | 100 |
| $s = 2$   | 69 | 79 | 84 | 98 |

We naturally observe that increasing the number of observations improves the performance of the procedure, with almost perfect results when $n = 5000$. However, the results obtained with small sample sizes ($n = 500$) are already encouraging (more than 65% of success).

In the case where the true parameter $s$ does not belong to the grid, we observed that the procedure recovers the value of the grid which is closest to $s$.

## 5. Proofs

We use $C$ to denote an absolute constant whose values may change along the lines.

*Proof of Proposition 1.* We fix the true value $s = s_k$ in the grid. Recall that the size of the grid is at most given by the step $d_n = c(\log n)^{-1}$. We want to control

$$\mathbb{P}_{f,s_k}(\hat{s}_n \neq s_k) = \mathbb{P}_{f,s_k}(\hat{s}_n > s_k) + \mathbb{P}_{f,s_k}(\hat{s}_n < s_k).$$

The overestimation case, namely $\hat{s}_n > s_k$, is the simplest to deal with. By definition of $\hat{s}_n$, we have,

$$\mathbb{P}_{f,s_k}(\hat{s}_n > s_k) \leq \mathbb{P}_{f,s_k}\left(|\hat{\Phi}_n^p(u_{n,k})| \leq \frac{1}{2}\{q_{\beta'}\Phi^{[k]} + \Phi^{[k+1]}\}(u_{n,k})\right)$$
$$+ \mathbb{P}_{f,s_k}\left(|\hat{\Phi}_n^p(u_{n,k})| \geq \frac{1}{2}\{q_{\beta'}\Phi^{[k-1]} + \Phi^{[k]}\}(u_{n,k})\right). \quad (14)$$

Considering the first term in the right hand side of the previous inequality, and using that $|\Phi^p(u_{n,k})| \geq q_{\beta'}(u_{n,k})\Phi^{[k]}(u_{n,k})$, we can write

$$\mathbb{P}_{f,s_k}\left(|\hat{\Phi}_n^p(u_{n,k})| \leq \frac{1}{2}\{q_{\beta'}\Phi^{[k]} + \Phi^{[k+1]}\}(u_{n,k})\right)$$

$$\leq \quad \mathbb{P}_{f,s_k}\left(|\hat{\Phi}_n^p(u_{n,k}) - \Phi^p(u_{n,k})| \geq |\Phi^p(u_{n,k})| - \frac{1}{2}\left\{q_{\beta'}\Phi^{[k]} + \Phi^{[k+1]}\right\}(u_{n,k})\right)$$

$$\leq \quad \mathbb{P}_{f,s_k}\left(|\hat{\Phi}_n^p(u_{n,k}) - \Phi^p(u_{n,k})| \geq \frac{1}{2}\left\{q_{\beta'}(u_{n,k})\Phi^{[k]}(u_{n,k}) - \Phi^{[k+1]}(u_{n,k})\right\}\right).$$

We will often use the following lemma.

**Lemma 1.** *For any $j \in \{1, \ldots, N-1\}$ and $k \in \{2, \ldots, N\}$, we have*

$$\frac{1}{2}\left\{q_{\beta'}\Phi^{[j]} - \Phi^{[j+1]}\right\}(u_{n,j}) \quad \geq \quad \frac{1}{2}q_{\beta'}(u_{n,j})\Phi^{[j]}(u_{n,j})(1 + o(1)),$$

$$and \quad \frac{1}{2}\left\{q_{\beta'}\Phi^{[k-1]} - \Phi^{[k]}\right\}(u_{n,k}) \quad \geq \quad \frac{1}{2}\Phi^{[k]}(u_{n,k})(1 + o(1)).$$

*Proof of Lemma 1.* By using both that $s_{j+1} - s_j \geq d_n$ and $d_n \log(u_{n,j}) \to 0$, we have

$$\frac{1}{2}\left\{q_{\beta'}\Phi^{[j]} - \Phi^{[j+1]}\right\}(u_{n,j})$$

$$= \quad \frac{1}{2}\Phi^{[j]}(u_{n,j})q_{\beta'}(u_{n,j})\left\{1 - A^{-1}u_{n,j}^{\beta'}\exp(u_{n,j}^{s_j} - u_{n,j}^{s_{j+1}})\right\}$$

$$\geq \quad \frac{1}{2}\Phi^{[j]}(u_{n,j})q_{\beta'}(u_{n,j})\left\{1 - A^{-1}u_{n,j}^{\beta'}\exp\left[u_{n,j}^{s_j}(1 - \exp(d_n\log u_{n,j}))\right]\right\}$$

$$= \quad \frac{1}{2}\Phi^{[j]}(u_{n,j})q_{\beta'}(u_{n,j})\left\{1 - A^{-1}u_{n,j}^{\beta'}\exp\left[-u_{n,j}^{s_j}d_n\log u_{n,j}(1 + o(1))\right]\right\}$$

$$= \quad \frac{1}{2}\Phi^{[j]}(u_{n,j})q_{\beta'}(u_{n,j})\left\{1 - A^{-1}\exp\left[(-c/2 + \beta')\log(u_{n,j})(1 + o(1))\right]\right\}$$

$$= \quad \frac{1}{2}\Phi^{[j]}(u_{n,j})q_{\beta'}(u_{n,j})(1 + o(1)),$$

as soon as $-c/2 + \beta' < 0$, i.e. $c > 2\beta'$. Similarly, we have

$$\frac{1}{2}\left\{q_{\beta'}\Phi^{[j-1]} - \Phi^{[j]}\right\}(u_{n,j})$$

$$= \quad \frac{1}{2}\Phi^{[j]}(u_{n,j})\left\{q_{\beta'}(u_{n,j})\exp(u_{n,j}^{s_j} - u_{n,j}^{s_{j-1}}) - 1\right\}$$

$$\geq \quad \frac{1}{2}\Phi^{[j]}(u_{n,j})\left\{A\exp(u_{n,j}^{s_j}d_n\log u_{n,j} - \beta'\log(u_{n,j})) - 1\right\}$$

$$= \quad \frac{1}{2}\Phi^{[j]}(u_{n,j})\left\{A\exp\left((c/2 - \beta')\log(u_{n,j})(1 + o(1))\right) - 1\right\}$$

$$\geq \quad \frac{1}{2}\Phi^{[j]}(u_{n,j})(1 + o(1)),$$

as soon as $c > 2\beta'$. This ends the proof of the lemma. $\qquad\square$

Using the first result of this lemma, combined with Hoeffding inequality, we obtain

$$\mathbb{P}_{f,s_k}\left(|\hat{\Phi}_n^p(u_{n,k})| \leq \frac{1}{2}\{q_{\beta'}\Phi^{[k]} + \Phi^{[k+1]}\}(u_{n,k})\right)$$

$$\leq \quad \mathbb{P}_{f,s_k}\left(|\hat{\Phi}_n^p(u_{n,k}) - \Phi^p(u_{n,k})| \geq \frac{1}{2}\Phi^{[k]}(u_{n,k})q_{\beta'}(u_{n,k})(1 + o(1))\right)$$

$$\leq \quad \exp\left(-\frac{n}{4}q_{\beta'}^2(u_{n,k})\exp(-2u_{n,k}^{s_k})(1 + o(1))\right)$$

$$\leq \quad \exp\left(-\frac{A^2}{4}2^{2\beta'/s_k}(\log n)^{\frac{2(\delta-\beta')}{s_k}}(1 + o(1)))\right).$$

Similarly, by using the bound $|\Phi^p(u_{n,k})| \leq \Phi^{[k]}(u_{n,k})$, and the second result of Lemma 1, the second term in the right hand side of (14) satisfies

$$\mathbb{P}_{f,s_k}\left(|\hat{\Phi}_n^p(u_{n,k})| \geq \frac{1}{2}\{q_{\beta'}\Phi^{[k-1]} + \Phi^{[k]}\}(u_{n,k})\right)$$

$$\leq \mathbb{P}_{f,s_k}\left(|\hat{\Phi}_n^p(u_{n,k}) - \Phi^p(u_{n,k})| \geq \frac{1}{2}\left\{q_{\beta'}\Phi^{[k-1]} - \Phi^{[k]}\right\}(u_{n,k})\right)$$

$$\leq \mathbb{P}_{f,s_k}\left(|\hat{\Phi}_n^p(u_{n,k}) - \Phi^p(u_{n,k})| \geq \frac{1}{2}\Phi^{[k]}(u_{n,k})(1 + o(1))\right).$$

Hoeffding inequality leads to

$$\mathbb{P}_{f,s_k}\left(|\hat{\Phi}_n^p(u_{n,k})| \geq \frac{1}{2}\{q_{\beta'}\Phi^{[k-1]} + \Phi^{[k]}\}(u_{n,k})\right)$$

$$\leq \exp\left(-\frac{1}{4}(\log n)^{2\delta/s_k}(1 + o(1))\right).$$

Finally, the overestimation probability (14) is bounded by

$$\mathbb{P}_{f,s_k}(\hat{s}_n > s_k) \leq \exp\left(-\frac{A^2}{4}2^{2\beta'/s_k}(\log n)^{\frac{2(\delta-\beta')}{s_k}}(1 + o(1))\right).$$

Let us now consider the probability of underestimation. The case $\hat{s}_n = s_1$ has to be dealt with separately as it may occur from emptyness of the set $\hat{S}_n$. By using the definition of $\hat{s}_n$, we have

$$\mathbb{P}_{f,s_k}(s_1 < \hat{s}_n < s_k)$$

$$\leq \mathbb{P}_{f,s_k}\left(\underset{2\leq j<k}{\cup}\left\{|\hat{\Phi}_n^p(u_{n,j})| \geq \frac{1}{2}\left\{q_{\beta'}\Phi^{[j]} + \Phi^{[j+1]}\right\}(u_{n,j})\right\}\right)$$

$$\leq \sum_{j=2}^{k-1}\mathbb{P}_{f,s_k}\left(|\hat{\Phi}_n^p(u_{n,j}) - \Phi^p(u_{n,j})| \geq \frac{1}{2}\left\{q_{\beta'}\Phi^{[j]} + \Phi^{[j+1]}\right\}(u_{n,j}) - |\Phi^p(u_{n,j})|\right).$$

As $|\Phi^p(u_{n,j})| \leq \Phi^g(u_{n,j}) = \Phi^{[k]}(u_{n,j}) \leq \Phi^{[j+1]}(u_{n,j})$, we get

$$\mathbb{P}_{f,s_k}(s_1 < \hat{s}_n < s_k)$$
$$\leq \sum_{j=2}^{k-1} \mathbb{P}_{f,s_k}\left(|\hat{\Phi}_n^p(u_{n,j})-\Phi^p(u_{n,j})| \geq \frac{1}{2}\left\{q_{\beta'}\Phi^{[j]}+\Phi^{[j+1]}\right\}(u_{n,j})-\Phi^{[j+1]}(u_{n,j})\right).$$

Now, using Lemma 1 again and the Hoeffding inequality

$$\mathbb{P}_{f,s_k}(s_1 < \hat{s}_n < s_k)$$
$$\leq \sum_{j=2}^{k-1} \mathbb{P}_{f,s_k}\left(|\hat{\Phi}_n^p(u_{n,j})-\Phi^p(u_{n,j})| \geq \frac{1}{2}\Phi^{[j]}(u_{n,j})q_{\beta'}(u_{n,j})(1+o(1))\right)$$
$$\leq \sum_{j=2}^{k-1} \exp\left(-\frac{n}{4}q_{\beta'}^2(u_{n,j})\exp(-2u_{n,j}^{s_j})\right)$$
$$\leq N\exp\left(-\frac{A^2}{4}2^{2\beta'/s_j}(\log n)^{\frac{2(\delta-\beta')}{s_j}}(1+o(1))\right)$$
$$\leq \exp\left(-\frac{A^2}{4}2^{2\beta'/\bar{s}}(\log n)^{\frac{2(\delta-\beta')}{\bar{s}}}(1+o(1))\right),$$

as $s_j < \bar{s}$ and $N = O(\log n)$.

The case $\hat{s}_n = s_1$ can now be easily handled. Indeed, let us denote by $\mathcal{E}_j$ the event

$$\mathcal{E}_j = \left\{|\hat{\Phi}_n^p(u_{n,1})| \geq 1/2\left(q_{\beta'}\Phi^{[j]}+\Phi^{[j+1]}\right)(u_{n,1})\right\}.$$

Now, if $\hat{s}_n = s_1$, then either the event $\mathcal{E}_1$ happens, or all of the $\mathcal{E}_j$s don't and thus in particular, $\mathcal{E}_k$ is not satisfied. Thus,

$$\mathbb{P}_{f,s_k}(\hat{s}_n = s_1) \leq \mathbb{P}_{f,s_k}(\mathcal{E}_1 \cup \mathcal{E}_k^c)$$

The probability of $\mathcal{E}_k^c$ has already been controlled (overestimation probability). Let us consider the probability of the first event. As previously seen, using Lemma 1 and Hoeffding inequality,

$$\mathbb{P}_{f,s_k}(\mathcal{E}_1)$$
$$\leq \mathbb{P}_{f,s_k}\left(|\hat{\Phi}_n^p(u_{n,1})-\Phi^p(u_{n,1})| \geq \frac{1}{2}\left\{q_{\beta'}\Phi^{[1]}+\Phi^{[2]}\right\}(u_{n,1})-|\Phi^p(u_{n,1})|\right)$$
$$\leq \mathbb{P}_{f,s_k}\left(|\hat{\Phi}_n^p(u_{n,1})-\Phi^p(u_{n,1})| \geq \frac{1}{2}\Phi^{[1]}(u_{n,1})q_{\beta'}(u_{n,1})(1+o(1))\right)$$
$$\leq \exp\left(-\frac{n}{4}q_{\beta'}^2(u_{n,1})\exp(-2u_n^{s_1})\right)$$
$$\leq \exp\left(-\frac{A^2}{4}2^{2\beta'/s_1}(\log n)^{\frac{2(\delta-\beta')}{s_1}}(1+o(1))\right).$$

Thus, the probability of underestimation is bounded by

$$\mathbb{P}_{f,s_k}(\hat{s}_n < s_k) \leq \exp\left(-\frac{A^2}{4}2^{2\beta'/\bar{s}}(\log n)^{2(\delta-\beta')/\bar{s}}(1+o(1))\right)$$

and gathering the results concerning overestimation and underestimation, we get

$$\mathbb{P}_{f,s_k}(\hat{s}_n \neq s_k) \leq \exp\left(-\frac{A^2}{4} 2^{2\beta'/\bar{s}}(\log n)^{2(\delta-\beta')/\bar{s}}(1 + o(1))\right)$$

$\Box$

*Proof of Corollary 1.* Let the true value of the parameter be some fixed point $s_k$ on the grid. We introduce respectively, $h_n$, the non-random version of the bandwidth $\hat{h}_n$ and $K_n$ the non-random version of the kernel $\hat{K}_n$ both constructed with true self-similarity index $s_k$. The Fourier transform $\Phi^{K_n}$ of $K_n$ thus satisfies

$$\Phi^{K_n}(u) = \exp\left\{\left(\frac{|u|}{\hat{h}_n}\right)^{s_k}\right\} 1_{|u| \leq 1},$$

$$\text{where } h_n = \left(\frac{\log n}{2} - \frac{\bar{\beta} - s_k + 1/2}{s_k} \log\log n\right)^{-1/s_k}.$$

We also introduce the corresponding (classical) estimator

$$f_n(x) = \frac{1}{nh_n} \sum_{i=1}^{n} K_n\left(\frac{x - Y_i}{h_n}\right),$$

which corresponds to the case of entirely known noise distribution. Note that obviously, $s_k, K_n$ and $h_n$ are unknown to the statistician. These objects are used only as tools to assess the convergence of the procedure. Now, remark that we have

$$\mathbb{E}_{f,s_k}[|\hat{f}_n(x) - f(x)|^2] = \mathbb{E}_{f,s_k}[|f_n(x) - f(x)|^2 1_{\hat{s}_n = s_k}]$$
$$+ \mathbb{E}_{f,s_k}[|\hat{f}_n(x) - f(x)|^2 1_{\hat{s}_n \neq s_k}] = T_1 + T_2,$$

say. Let us focus on the first term

$$T_1 \leq \mathbb{E}_{f,s_k}[|f_n(x) - f(x)|^2] = \{\mathbb{E}_{f,s_k}[f_n(x)] - f(x)\}^2 + \mathbb{V}\text{ar}_{f,s_k}\{f_n(x)\},$$

introducing the bias and the variance of the estimator $f_n(x)$. By using classical results on this estimator, we have

$$T_1 \leq O(h_n^{2\beta-1}) + O\left(\frac{h_n^{2(s_k-1)}\exp(2/h_n^{s_k})}{n}\right).$$

Now, we prove that the second term $T_2$ is negligible in front of the main term $T_1$, by using Proposition 1 and uniform bounds on $|\hat{f}_n(x)|$ and $|f(x)|$. First,

$$|\hat{f}_n(x)| \leq \int e^{|t|^{\bar{s}}} 1_{|t| \leq 1/\hat{h}_n} dt = O(\hat{h}_n^{\bar{s}-1} \exp\{1/\hat{h}_n^{\bar{s}}\})$$
$$\leq O(1)(\log n)^{(1-\bar{s})/\underline{s}} \exp\{(\log n)^{\bar{s}/\underline{s}}\}$$

and also

$$|f(x)| \leq \int |\Phi^f(t)| dt = O(\int (1 + |t|^{2\beta})^{-1} dt) = O(1),$$

leading to

$$T_2 = O((\log n)^{2(1-\bar{s})/\underline{s}} \exp\{2(\log n)^{\bar{s}/\underline{s}}\}) \mathbb{P}_{f,s_k}(\hat{s}_n \neq s_k)$$

$$= O\left((\log n)^{2(1-\bar{s})/\underline{s}} \exp\left(2(\log n)^{\bar{s}/\underline{s}} - \frac{A^2}{4} 2^{2\beta'/\bar{s}} (\log n)^{2(\delta - \beta')/\bar{s}} (1 + o(1))\right)\right).$$

As soon as we choose $2(\delta - \beta')/\bar{s} > \bar{s}/\underline{s}$, this second term $T_2$ will be negligible in front of $T_1$. In conclusion,

$$\mathbb{E}_{f,s_k}[|\hat{f}_n(x) - f(x)|^2] = O(h_n^{2\beta-1}) + O\left(h_n^{2(s_k-1)} \frac{\exp(2/h_n^{s_k})}{n}\right)$$
$$= O((\log n)^{-(2\beta-1)/s_k}).$$

$\square$

*Proof of Corollary 2.* We keep on with the same notations as in the proof of Corollary 1 and denote by $I$ the functional $\int f^2$ and by $T_n$ the estimator using deterministic parameters $s_k, K_n$ and $h_n$. In the same way as in the proof of Corollary 1, we write

$$\mathbb{E}_{f,s_k}[|\hat{T}_n - I|^2] \leq \mathbb{E}_{f,s_k}[|T_n - I|^2] + \mathbb{E}_{f,s_k}[|\hat{T}_n - I|^2 1_{\{\hat{s}_n \neq s_k\}}]. \tag{15}$$

Let us first focus on the first term appearing in the right hand side of (15). We split it into the square of a bias term plus a variance term. The bias is bounded by

$$|\mathbb{E}_{f,s_k} T_n - I| \leq O((\log n)^{-2\beta/s_k}).$$

Concerning the variance term, we easily get

$$\mathbb{V}\text{ar}_{f,s_k}(T_n) \leq \frac{C_1}{n^2} h_n^{s_k-1} \exp(4/h_n^{s_k}) + \frac{C_2}{n} h_n^{2\beta+s_k-1} \exp(2/h_n^{s_k}),$$

where $C_1$ and $C_2$ are positive constants (we refer to [1], Theorem 4 for more details). Using the form of the bandwidth $h_n$, we have

$$\mathbb{E}_{f,s_k} |T_n - I|^2 = O\left(\frac{\log n}{2}\right)^{-4\beta/s_k}.$$

Let us now focus on the second term appearing in the right hand side of (15). Denoting by $h_0 = (\log n/2)^{-1/\underline{s}}$, we have

$$|\hat{T}_n| \leq \frac{1}{2\pi} \int_{|u| \leq 1/h_0} \exp(2|u|^{\bar{s}}) du = O(h_0^{\bar{s}-1} \exp(2/h_0^{\bar{s}})).$$

Moreover,

$$I = \|f\|_2^2 = \frac{1}{2\pi} \|\Phi^f\|_2^2$$

This leads to

$$\mathbb{E}_{f,s_k}[|\hat{T}_n - I|^2 \mathbf{1}_{\{\hat{s}_n \neq s_k\}}]$$

$$\leq C \left(\frac{\log n}{2}\right)^{(1-\bar{s})/\underline{s}} \exp\left\{2\left(\frac{\log n}{2}\right)^{\bar{s}/\underline{s}}\right\} \mathbb{P}_{f,s_k}(\hat{s}_n \neq s_k)$$

$$\leq C \left(\frac{\log n}{2}\right)^{(1-\bar{s})/\underline{s}} \exp\left\{2\left(\frac{\log n}{2}\right)^{\bar{s}/\underline{s}}\right\}$$

$$\times \exp\left(-\frac{A^2}{4} 2^{2\beta'/\bar{s}} (\log n)^{2(\delta-\beta')/\bar{s}} (1 + o(1))\right),$$

and this term is negligible in front of the first term appearing in the right hand side of (15) as soon as $2(\delta - \beta')/\bar{s} > \bar{s}/\underline{s}$. This leads to the result. $\qquad\square$

*Proof of Corollary 3.* We use the same notations as in the proof of Corollaries 1 and 2. Moreover, $T_n^0$ is the test statistic constructed with the deterministic kernel $K_n$ and the deterministic bandwidth $h_n$; and $t_n^2$ is the threshold defined with the true parameter value $s_k$ for the self-similarity index. The first type error of the test is controlled by

$$\mathbb{P}_{f_0,s_k}(\Delta_n^\star = 1) = \mathbb{P}_{f_0,s_k}(|\hat{T}_n^0|\hat{t}_n^{-2} > \mathcal{C}^\star) \leq \mathbb{P}_{f_0,s_k}(\hat{s}_n \neq s_k) + \mathbb{P}_{f_0,s_k}(|T_n^0|t_n^{-2} > \mathcal{C}^\star).$$

The first term on the right hand side of this inequality converges to zero according to Proposition 1. Moreover, Theorem 4 in [1], shows that

$$\mathbb{V}\mathrm{ar}_{f_0,s_k}(T_n^0) \leq O(1)\frac{h_n^{s_k-1}}{n^2}\exp(4/h_n^{s_k}) + O(1)\frac{h_n^{2\bar{\beta}+s_k-1}}{n}\exp(2/h_n^{s_k}).$$

Finally, we get

$$\mathbb{P}_{f_0,s_k}(|T_n^0|t_n^{-2} > \mathcal{C}^\star) \leq \frac{1}{(\mathcal{C}^\star)^2 t_n^4}$$

$$\times \left\{O(h_n^{4\bar{\beta}}) + O(1)\frac{h_n^{s_k-1}}{n^2}\exp(4/h_n^{s_k}) + O(1)\frac{h_n^{2\bar{\beta}+s_k-1}}{n}\exp(2/h_n^{s_k})\right\},$$

which is actually $O(1)/\mathcal{C}^\star$. Choosing $\mathcal{C}^\star$ large enough achieves the control of the first error term.

We now turn to the second type error term. Under hypothesis $H_1(\mathcal{C}, \Psi_n)$, there exists some $\beta$ such that $f$ belongs to $\mathcal{S}(\beta, L)$ and $\|f - f_0\|_2^2 \geq \mathcal{C}\psi_{n,\beta}$. We write

$$\mathbb{P}_{f,s_k}(\Delta_n^\star = 0) = \mathbb{P}_{f,s_k}(|\hat{T}_n^0|\hat{t}_n^{-2} \leq \mathcal{C}^\star) \leq \mathbb{P}_{f,s_k}(\hat{s}_n \neq s_k) + \mathbb{P}_{f,s_k}(|T_n^0|t_n^{-2} \leq \mathcal{C}^\star).$$

As already seen, the first term in the right hand side of this inequality converges to zero, so we only deal with the second one. We define $B_{f,s_k}(T_n^0) = \mathbb{E}_{f,s_k}T_n^0 - \|f - f_0\|_2^2$. Thus

$$\mathbb{P}_{f,s_k}(|T_n^0|t_n^{-2} \leq \mathcal{C}^\star) \leq \mathbb{P}_{f,s_k}(|T_n^0 - \mathbb{E}_{f,s_k}T_n^0| \geq \|f - f_0\|_2^2 - \mathcal{C}^\star t_n^2 + B_{f,s_k}(T_n^0))$$

$$\leq \frac{\mathbb{V}\mathrm{ar}_{f,s_k}(T_n^0)}{(\|f - f_0\|_2^2 - \mathcal{C}^\star t_n^2 + B_{f,s_k}(T_n^0))^2}. \qquad (16)$$

According to [1], we have

$$B_{f,s_k}(T_n^0) \leq C_1 h_n^{2\beta}$$

where $C_1 > 0$ is a constant depending only on $L$ and on the noise distribution. Under hypothesis $H_1(\mathcal{C}, \Psi_n)$, we also have $\|f - f_0\|_2^2 \geq \mathcal{C}\psi_{n,\beta}^2$. Thus,

$$\|f - f_0\|_2^2 - \mathcal{C}^\star t_n^2 + B_{f,s_k}(T_n^0)$$
$$\geq \mathcal{C}\left(\frac{\log n}{2}\right)^{-2\beta/s_k} - \mathcal{C}^\star\left(\frac{\log n}{2}\right)^{-2\tilde{\beta}/s_k} - C_1\left(\frac{\log n}{2}\right)^{-2\beta/s_k}$$
$$\geq a\left(\frac{\log n}{2}\right)^{-2\beta/s_k}.$$

where $a = \mathcal{C} - \mathcal{C}^\star - C_1$ is positive whenever $\mathcal{C} > \mathcal{C}^0 := \mathcal{C}^\star - C_1$. Returning to (16), we get

$$\mathbb{P}_{f,s_k}(|T_n^0|t_n^{-2}) \leq \frac{\psi_{n,\beta}^4}{a^2}\mathbb{V}\mathrm{ar}_{f,s_k}(T_n^0).$$

Computation of the variance follows the same lines as under hypothesis $H_0$. We obtain

$$\mathbb{V}\mathrm{ar}_{f,s_k}(T_n^0) \leq O(1)\frac{h_n^{s_k-1}}{n}\exp(2/h_n^{s_k})\left(h_n^{2\beta} + \frac{\exp(2/h_n^{s_k})}{n}\right).$$

The choice of the bandwidth ensures that the second type error term converges to zero.    □

## References

[1] Butucea, C. (2007). Goodness-of-fit testing and quadratic functional estimation from indirect observations. *Ann. Statist.* **35**, 5, 1907–1930.

[2] Butucea, C. and Comte, F. (2008). Adaptive estimation of linear functionals in the convolution model and applications. *Bernoulli*. To appear.

[3] Butucea, C. and Matias, C. (2005). Minimax estimation of the noise level and of the deconvolution density in a semiparametric convolution model. *Bernoulli* **11**, 2, 309–340. MR2132729

[4] Butucea, C., Matias, C., and Pouet, C. (2008). Adaptive goodness-of-fit testing from indirect observations. *Annales de l'Institut Poincaré, Probabilités et Statistiques*. To appear.

[5] Butucea, C. and Tsybakov, A. B. (2007). Sharp optimality for density deconvolution with dominating bias. I. *Theory Probab. Appl* **52**, 1, 111–128.

[6] Efromovich, S. (1997). Density estimation for the case of supersmooth measurement error. *J. Amer. Statist. Assoc.* **92**, 438, 526–535. MR1467846

[7] Fan, J. (1991). On the optimal rates of convergence for nonparametric deconvolution problems. *Ann. Statist.* **19**, 3, 1257–1272.

[8] Holzmann, H., Bissantz, N., and Munk, A. (2007). Density testing in a contaminated sample. *J. Multivariate Analysis* **98**, 1, 57–75. MR2292917

[9] Ingster, Y. and Suslina, I. (2003). *Nonparametric goodness-of-fit testing under Gaussian models.* Lecture Notes in Statistics. 169. Springer, New York. MR1991446

[10] Matias, C. (2002). Semiparametric deconvolution with unknown noise variance. *ESAIM, Probab. Stat. 6*, 271–292. (electronic). MR1943151

[11] Meister, A. (2004). Deconvolution density estimation with a testing procedure for the error distribution. Tech. rep., Universität Stuttgart, 2004-004. http://www.mathematik.uni-stuttgart.de/preprints/listen/listen.php.

[12] Meister, A. (2006). Density estimation with normal measurement error with unknown variance. *Statist. Sinica* **16**, 1, 195–211.

[13] Neumann, M. H. (1997). On the effect of estimating the error density in nonparametric deconvolution. *J. Nonparametr. Statist.* **7**, 4, 307–330. MR1460203

[14] Weron, R. (1996). On the Chambers-Mallows-Stuck method for simulating skewed stable random variables. *Statist. Probab. Lett.* **28**, 2, 165–171. MR1394670 (97c:60037)